\documentclass[12pt]{article}
\usepackage[log-declarations=false]{xparse}
\usepackage{amsmath,amsxtra,amssymb,latexsym, amscd,amsthm}
\usepackage{verbatim}
\usepackage[colorlinks=true,citecolor=black,linkcolor=black,urlcolor=blue]{hyperref}
\usepackage{enumitem}

\usepackage[abbrev,logical-quotes,msc-links]{amsrefs}

\newtheorem{theorem}{Theorem}

\newtheorem{lemma}[theorem]{Lemma}

\newtheorem{definition}[]{Definition}
\newtheorem{conjecture}[]{Conjecture}

\newcommand{\bigO}{O}	
\newcommand{\tricross}{C}	

\theoremstyle{plain}

\newcommand{\stack}[6]{$$ \mbox{ #1 } \mbox{ #2 } \left\{
	   \begin{array}{ll}
	   \mbox{ #3 } & \mbox{ #4 } \\ \\
	   \mbox{ #5 } & \mbox{ #6 }
	   \end{array} \right. $$ }

\newcommand{\nmsystem}{$(n,m)$-structure}		

 {\begin{list}{}%
         {\setlength{\leftmargin}{#1}}%
         \item[]%
 }
 {\end{list}}

\newcounter{startSubdivision}
\newcounter{subdividedLemma}
\renewcommand*{\thesubdividedLemma}{\arabic{theorem}} 

\newenvironment{subdividedLemma}{%
	\refstepcounter{theorem}%
	\paragraph{Lemma~\thesubdividedLemma.}%
	\enumerate[start=\value{startSubdivision}]
}{%
	\setcounter{startSubdivision}{\value{enumi}}
	\refstepcounter{startSubdivision}
	\endenumerate
}

\title{An asymptotic resolution \\of a conjecture of Szemer\'edi and Petruska}

\author{Andr\'{e} E. K\'{e}zdy \\
	\small Department of Mathematics\\[-0.8ex]
	\small University of Louisville\\[-0.8ex]
	\small Louisville, Kentucky, U.S.A.\\[-0.8ex]
	\small\tt kezdy@louisville.edu\\
	\\
	Jen\H{o} Lehel \\
\small Department of Mathematics\\[-0.8ex]
	\small University of Louisville\\[-0.8ex]
	\small Louisville, Kentucky, U.S.A.\\
                \small{and}\\          
           \small  Alfr\'ed R\'enyi Institute of Mathematics\\[-0.8ex] 
            \small  Budapest, Hungary\\[-0.8ex]
              	\small\tt lehelj@renyi.hu}  
              
\begin{document}
\maketitle

\begin{abstract}
Consider a $3$-uniform hypergraph of order $n$ with
clique number $k$ such that the intersection of all its $k$-cliques is empty. Szemer\'edi and Petruska  
proved  $n\leq 8m^2+3m$, for fixed $m=n-k$, and they 
 conjectured the sharp bound $n \leq {m+2 \choose 2}$.  
This problem is known to be equivalent
to determining the maximum order of a $\tau$-critical $3$-uniform hypergraph with transversal number $m$ 
(details may also be found in a companion paper  \cite{KL}).

The best known bound, $n\leq \frac{3}{4}m^2+m+1$, was obtained by Tuza using the machinery of $\tau$-critical hypergraphs. 
Here we propose an alternative approach, a combination  of the iterative  decomposition process introduced by Szemer\'edi and Petruska 
with the skew version of Bollob\'as's theorem on set pair systems. The new approach improves the bound to
$n\leq {m+2 \choose 2} + \bigO(m^{{5}/{3}})$, resolving the conjecture asymptotically.
\end{abstract}

\section{Introduction}
\label{Introduction}

Let ${\cal N} = \{N_1,\ldots,N_\ell\}$ be a collection of $k$-subsets of $ \{1,\ldots,n\}$.
Set $V = \bigcup_{i=1}^\ell N_i$.
Assume that $n = |V|$, $\ell \geq 2$, and $k \geq 3$.  Set $m = n - k$; that is, $\left|\overline{N_i}\right| =|V \setminus N_i|= m$.
We further assume that ${\cal N}$ satisfies the following two properties:
\begin{quote}
\begin{enumerate}[label=({\it\roman*\hspace{.1em}})]  
\item \label{MinimalNonIntersectingProperty} $\bigcap_{i=1}^\ell N_i = \varnothing$, but $\bigcap_{j\neq i} N_j \neq \varnothing$ for all $i=1,\ldots,\ell$. 
\item \label{LargestCliqueSizeProperty} For any $S \subseteq V$ such that $|S| \geq k+1$, there exists a subset $T \subseteq S$ such that $|T| = 3$ and 
$T \not\subseteq N_i$, for all $i=1,\ldots,\ell$.
\end{enumerate}
\end{quote}
\noindent 
Notice that property \ref{MinimalNonIntersectingProperty} means that ${\cal N}$ is minimal with respect to having empty intersection; and 
property \ref{LargestCliqueSizeProperty} may be interpreted as stating that, in the $3$-uniform hypergraph induced by all $3$-subsets of the $k$-sets in ${\cal N}$, the largest cliques have $k$ vertices.

We shall refer to a collection ${\cal N}$ satisfying \ref{MinimalNonIntersectingProperty} and \ref{LargestCliqueSizeProperty} as an \emph{\nmsystem{}}. Szemer\'edi and Petruska \cite{SzemerediPetruska} conjectured the following:

\begin{conjecture}\label{SPConjecture} Any \nmsystem{} satisfies $n \leq {m+2 \choose 2}$.
\end{conjecture}

Szemer\'edi and Petruska give the following construction to show that Conjecture \ref{SPConjecture}, if true, would be sharp.
For fixed integers $k,m\geq 3$ such that $k-1={m+1\choose 2}$, begin with disjoint sets $A$ and $G$ satisfying 
$|A|=k-1$ and  $|G|=m+1$. Let $A=\{a_1,\dots,a_{k-1}\}$
and let $\{p_1,\dots, p_{k-1}\}$  be the set of pairs ($2$-subsets) of $G$.
Define $N_i=(A\setminus\{a_i\})\cup p_i$, for $i=1,\dots,  k-1$. It is easy to check that the collection $\mathcal{N}=\{N_1,\dots,N_{k-1}\}$ satisfies properties (i) and (ii). 
 In particular,  the {\nmsystem{}}
$\mathcal{N}$ induces a $3$-uniform  hypergraph $H$ on vertex set 
$A\cup G$ with all $3$-subsets included by some $k$-set of $\cal{N}$ forming the edge set; 
furthermore, 
the order of $H$ is $n=|A|+|G|=k+m={m+2\choose 2}$ and $\omega(H)=k$. 

For $m=2,3$, there are other extremal \nmsystem{}s; however,
it has been conjectured (by us and others) that this construction is the unique extremal structure for $m\geq 4$.  For $m=2,3$ and $4$, Conjecture \ref{SPConjecture}  has been verified and all extremal structures have been characterized by Jobson et al. \cite{JobsonKezdyLehel}. 

The best known general bound, $n\leq \frac{3}{4}m^2+m+1$, was obtained by Tuza\footnote{ Personal communication ($2019$).} 
using the machinery of $\tau$-critical hypergraphs.  
Here we significantly refine and vindicate an alternative approach first proposed
by Jobson et al. \cite{JobsonKezdyPervenecki}. 
This approach develops the iterative decomposition process introduced by Szemer\'edi and Petruska and 
applies the skew version of Bollob\'as's theorem \cite{Bollobas}   on the size of cross-intersecting 
set pair 
systems. The result (Theorem \ref{maintheorem}) is an asymptotically tight upper bound:  
$$n\leq  {m+2 \choose 2} + 6 m^{5/3} + 3m^{4/3} + 9m - 3.$$

As noted by Gy\'arf\'as et al. 
\cite{GyarfasLehelTuza}, the Szemer\'edi and Petruska problem is equivalent
to determining the maximum order of a $\tau$-critical $3$-uniform hypergraph with transversal number $m$.  
More generally, they also determined that $\bigO(m^{r-1})$ is the correct order of magnitude for the
maximum order of a $\tau$-critical $r$-uniform hypergraph with transversal number $m$; 
the best known bounds were obtained by Tuza \cite{Tuza}.  

A companion paper \cite{KL} shows the details of the mentioned equivalence and presents further remarks on the origin of the Szemer\'edi and Petruska problem.  As an immediate corollary of our main theorem and this equivalence  we obtain a new and asymptotically tight bound for the order of a $\tau$-critical $3$-uniform hypergraph.

\begin{theorem}
If $H$ is a $3$-uniform $\tau$-critical hypergraph with $\tau(H)=t$, then $$|V(H)|\leq {t+2\choose 2}+\bigO(t^{5/3}).$$
\end{theorem}

Section \ref{decomposition} introduces notations and
recalls the process, introduced by Szemer\'edi and Petruska, to decompose an \nmsystem{} into stages. Consecutive stages can be viewed as survival times of fewer and fewer $k$-sets of the system.  
The basic concepts associated with an \nmsystem{} are the `kernels' belonging to the surviving $k$-sets at each time and the `private pairs' selected for the remaining $k$-subsets at each stage  
(the private pairs are pairs that belong to precisely one $k$-subset of the surviving subsystem at a stage of the decomposition). 
For example, in the conjectured extremal Szemer\'edi-Petruska construction described above, there is only one stage; the
set $A$ is the kernel of the only stage, and the $2$-subsets of $G$ are the private pairs of the $k$-sets. 
In Sections \ref{PairSelection}, \ref{advancement} and \ref{digraph} the iterative decomposition process used by Szemer\'edi and Petruska is extended considerably. 

In  addition to extending the iterative decomposition process our new approach 
applies the skew version of Bollob\'as's theorem \cite{Bollobas}  on the size of cross-intersecting 
set pair systems (Theorem \ref{SkewBollobas}):   
If $A_1,\ldots,A_h$ are $r$-element sets and $B_1,\ldots,B_h$ are $s$-element sets such that
$A_i\cap B_i=\varnothing$ for $i=1,\ldots,h$, and 
 $A_i\cap B_j\neq\varnothing $ whenever $1 \leq i < j \leq h$,
then $h \leq {r+s \choose r}$. 
Theorem \ref{SkewBollobas}  will be applied in Section  \ref{Skew}  with $r=2$ and $ s=m$,  
where the $2$-sets are  carefully selected private pairs, and  the  $m$-sets are derived iteratively from  the $k$-sets corresponding to all those private pairs. 

Section \ref{PairSelection} introduces a recursive procedure, based on the decomposition process,
to select special private pairs by `promoting' the initial pairs as needed.  
Section \ref{advancement} describes an additional processing of these special private pairs, known as `advancement',
which sacrifices or deactivates  a limited number of private pairs in order to guarantee that others may advance into a protected position
before choosing `free' private pairs in Sections \ref{digraph} and \ref{transversal}.

Sections \ref{digraph} and \ref{transversal} define a large subset  of free private pairs chosen 
from the special advanced pairs obtained in Section \ref{advancement}. A skew cross-intersecting $(2,m)$-system ultimately arises from this subset
of free private pairs and by appropriately adjusting the complements of the corresponding $k$-sets. This is done by
 using the recursive process (\ref{RecursiveDefinitionOfMs}) in Section \ref{Skew},  
where the proof of the main result, Theorem \ref{maintheorem}, 
is concluded.  

\section{The Decomposition Process} 
\label{decomposition}

We begin by giving definitions and recalling the process, introduced by Szemer\'edi and Petruska\footnote{ We have endeavored 
to use the same notation introduced by Szemer\'edi and Petruska.  Some important exceptions: we use $a_i^{(j)}$ for their $x_i^{(j)}$; also $k$ and $\ell$ here refer to their quantities $n$ and $k$, respectively.}, to decompose \nmsystem{}s.  
Much of this section is very similar to their presentation.
We assume $\ell\geq 4$ since Szemer\'edi and Petruska resolve the cases $\ell=2,3$.
Let ${\cal N} = \{N_1,\ldots,N_\ell\}$, be an \nmsystem{}.
Define a collection of objects in {\em stages}, which are also called {\em times}, starting with stage $0$.
Set $\ell_0 = \ell$, ${\cal N}^{(0)} = {\cal N}$ and $N_i^{(0)} = N_i$.  
For every $i=1,\ldots,\ell_0$, fix a choice of vertex $a_i^{(0)} \in \bigcap_{j\neq i} N_j$.  By definition, $a_i^{(0)} \neq a_j^{(0)}$, for $i \neq j$.  
The set $A^{(0)} = \left\{ a_1^{(0)},\ldots,a_{\ell_0}^{(0)} \right\}$ is called the {\em kernel} and  $a_1^{(0)},\ldots,a_{\ell_0}^{(0)}$ are called the {\em kernel vertices} at stage $0$.

Assume 
that  the collection of $k$-subsets ${\cal N}^{(j)} = \left\{N_1^{(j)},\ldots,N_{\ell_j}^{(j)}\right\}$  and  a corresponding kernel $A^{(j)} = \left\{ a_1^{(j)},\ldots,a_{\ell_j}^{(j)} \right\}$   
are defined  for  all  stages $0,1,\ldots,j$. Also assume that  $\ell_j\geq 4$ and the  sets in the 
`remainder' structure
$$
R^{(j)} = \left\{N_1^{(j)} \setminus \bigcup_{i=0}^j A^{(i)},\ldots,N_{\ell_j}^{(j)}\setminus \bigcup_{i=0}^j A^{(i)}\right\}
$$
have no common vertex.  We refer to $N_r^{(j+1)}=N_r^{(j)} \setminus \bigcup_{i=0}^j A^{(i)}$ as the {\em truncation} of $N_r^{(j)}$. 
We now explain the definition of 
${\cal N}^{(j+1)}$ and $A^{(j+1)}$.  

Because the truncations of the $N_r^{(j)}$'s in $R^{(j)}$ 
have no common vertex, there exist substructures of $R^{(j)}$ 
satisfying property \ref{MinimalNonIntersectingProperty}. 
Stop if $R^{(j)}$ contains such substructure(s) only with two or three sets.  Otherwise, let
$$
{\cal N}^{(j+1)} = \left\{N_1^{(j+1)},\ldots,N_{\ell_{j+1}}^{(j+1)}\right\}\subset R^{(j)} 
$$
be a substructure satisfying \ref{MinimalNonIntersectingProperty} and $\ell_{j+1} \geq 4$. 

For $i=1,\ldots,\ell_{j+1}$, fix a choice of vertex $a_i^{(j+1)} \in \bigcap\limits_{r=1 \atop r\neq i}^{\ell_{j+1}} N_r^{(j+1)} $ 
and let $A^{(j+1)} = \left\{ a_1^{(j+1)},\ldots,a_{\ell_{j+1}}^{(j+1)} \right\}$ be the kernel at time $j+1$. 
Observe that the sets in the remainder structure 
$$R^{(j+1)}=\left\{N_1^{(j+1)} \setminus \bigcup_{i=0}^{j+1} A^{(i)},\ldots,N_{\ell_{j+1}}^{(j+1)}\setminus \bigcup_{i=0}^{j+1} A^{(i)}\right\}$$
have empty intersection.

This process terminates at some time $t$, and defines  $\ell_{j}$, ${\cal N}^{(j)}$, and $A^{(j)}$ for $j=0,1, \ldots, t$.   The only substructures of the terminal remainder structure, $R^{(t)}$, that satisfy property (i) 
have $2$ or $3$ sets.  
Because ${\cal N} = \left\{N_1,\ldots,N_{\ell_0}\right\}$ is an arbitrary enumeration of ${\cal N}$, we may assume that
$$
{\cal N}^{(j)} = \left\{N_1,\ldots,N_{\ell_{j}}\right\}, \mbox{ for } j = 0,\ldots,t.
$$
That is, we linearly order the $k$-sets in ${\cal N}$ according to how long their truncations appear in the sequence of chosen substructures.   The longer that its truncations appear, the earlier a set appears in this linear ordering.

Define, for $i=1,\ldots,\ell_0$, the last {\em time} (or stage), denoted $t_i$,  that truncations of $N_i$ appear 
in a substructure of this decomposition process; equivalently, set
$$
t_i = \max\left\{j : N_i \in {\cal N}^{(j)}\right\}.
$$
By definition, $\ell= \ell_0\geq \cdots \geq \ell_{t} \geq 4$ and $t = t_1 \geq \cdots \geq t_{\ell_0} \geq 0$. 
Observe that $x\in\{1,\ldots,\ell_{j}\}$ implies $j \leq t_x$.

Define $A_j = \bigcup_{s=0}^{j} A^{(s)}$ and let $A=A_t$ be the set of all kernel vertices.  Let $G = V \setminus A$ denote the {\em garbage} vertices; that is, the vertices remaining after the  decomposition process terminates. 
The selected linear ordering of ${\cal N}$ together with the decomposition process induces a natural ordering of $V$:   
vertices within a kernel $A^{(s)}$ are ordered according to their subscripts ($a_i^{(s)} < a_j^{(s)}$ if and only if $i < j$),
kernels are ordered according to stage ($A^{(r)} < A^{(s)}$ if and only if $r < s$), and garbage vertices are ordered (linearly) last (in no particular order).

The next lemma begins a list of useful properties of the decomposition process.  
The results mentioned in Lemma \ref{FirstLemma} were proved in 
Lemma 5 and in Lemma 6 of \cite{SzemerediPetruska}.  We reprove
these results here for completeness.

\begin{subdividedLemma} \label{FirstLemma}  
	\item\label{tUpperBound} $t<m\leq |G|$
	\item\label{GUpperBound}  $|G| \leq 3(m-t-1)$.
\end{subdividedLemma}

\begin{proof} (\ref{tUpperBound}) Since $N_1\in {\cal N}^{(t)}$, we have  $\left|\overline{N_1}\cap A \right|=t+1$; and because  
    $|\overline{N_1}|=m$, it follows that $t<m$. Next suppose, on the contrary, that
$|G|\leq m-1$. Observe that $|A|=|V|-|G|\geq n-m+1=k+1$.  
 By property \ref{LargestCliqueSizeProperty}, there exists a set
 $T\subset A$, $|T|=3$, such that 
 $T\not\subseteq N_i$, for all $1\leq i\leq \ell$. 
Each element of $T$ is missed by at most one $N_i$, $1\leq i\leq \ell_t$.
It follows that because $\ell_t\geq 4$, one of these sets contains $T$, a contradiction. 
\vskip.5em
\noindent (\ref{GUpperBound}) 
When the decomposition process terminates  
any subsystem of
$$
R^{(t)}=\left\{ N_1\setminus A, \ldots, N_{\ell_t}\setminus A\right\}
$$
satisfying \ref{MinimalNonIntersectingProperty} has at most three sets; we may assume that these are $\{N_1\setminus A,N_2\setminus A,N_3\setminus A\}$ or $\{N_1\setminus A,N_2\setminus A\}$. Each set $N_i$, $i=1,2$ or $3$, `survives'  until stage $t$, hence  
$$r=\left|N_i \cap G\right| = \left|N_i \setminus  A\right|=k - (\ell_0-1) - \cdots - (\ell_t-1) = k - |A| + t + 1.$$
Thus, the subsystem consists of two or three 
$r$-sets that have empty intersection; so
the complements of these sets in $G$
must satisfy $\left|\overline{N_i}\cap G \right|\geq \frac{r}{2}$.  Consequently, 
$$|G| = \left|N_1 \cap G\right| + \left|\overline{N_1}\cap G \right|\geq r + \frac{r}{2}  = \frac{3}{2} (k - |A| + t + 1).$$

   By substituting $k=n-m$ and $n-|A|=|G|$, we obtain    
 $$|G|\geq   \frac{3}{2} (k - |A| + t + 1)
 =\frac{3}{2} (-m +|G|+ t + 1),$$
thus $|G| \leq 3(m-t-1)$ follows.
\end{proof}

\section{Selection of private pairs} 
\label{PairSelection}

For $j$, $0 \leq j \leq t_i$, a pair of vertices $p \subset N_i$ is {\em single-covered with respect to ${\cal N}^{(j)}$} 
if $N_i$ is the only set in ${\cal N}^{(j)}$ that contains $p$ as a subset. 
A pair that is contained in at least two sets in ${\cal N}^{(j)}$ is called a
{\em double-covered} pair (at time $j$). If $p \subset N_i$ is single-covered with respect to ${\cal N}^{(j)}$,
then it is called a {\em private pair for $N_i$ at time $j$}, or simply a {\em private pair} of $N_i$.
We also simply say that ``$p$ is private for $N_i^{(j)}$'' to mean ``$p$ is a private pair for $N_i$ at time $j$''.
Observe that, if a pair
is private for $N_i^{(j)}$, then it remains a private pair for $N_i$ until (and including) time $t_i$.
Similar terminology is used for private elements.  For example, a vertex 
$v$ is private for $N_i^{(j)}$ if $N_i$ is the only set in ${\cal N}^{(j)}$
that contains $v$.

The following lemma is a rephrasing of Lemma 7(a) proven 
by Szemer\'edi and Petruska \cite{SzemerediPetruska}.  
Recall the notation  $A_j=\bigcup_{s=0}^j A^{(s)}$.

\begin{subdividedLemma} \label{PairsInKernel} 
	\item\label{PairsInKernelAreDoubleCovered} For all $j=0,\ldots,t$, every pair $p\subset A_j$ is double-covered at time $j$.
\end{subdividedLemma}
\begin{proof}
    (\ref{PairsInKernelAreDoubleCovered}) 
    Suppose that $p =\{a_{i_1}^{({j_1})}, a_{i_2}^{({j_2})}\} \subset A_j$, for some $j \in \{0,\ldots,t\}$.
    Necessarily $j_1, j_2 \leq j$.
    Because $\ell_j \geq 4$, there exists $\{i_3,i_4\} \subseteq \{ 1,\ldots,\ell_j\} \setminus \{i_1,i_2\}$.
    By definition of the kernel, $p \subseteq N_{i_3}^{(j)} \cap N_{i_3}^{(j)}$ so $p$ is double-covered at time $j$.
\end{proof}

The next three lemmas introduce new and powerful tools; they are a crucial refinement of the proof of Szemer\'edi and Petruska's
Lemma 7(b) \cite{SzemerediPetruska}.  
In particular, Lemma \ref{ExistenceOfThreeCrosses} introduces the novel notion of a $3$-cross.

\begin{subdividedLemma} \label{ExistenceOfThreeCrosses} 
	\item[] Suppose that $N_i \in {\cal N}^{(j)}$, $Y \subseteq N_i$ and $|Y| \leq j$.
	\item\label{ThreeCrossExists} There exists a $3$-set $\tricross_i^{(j)}(Y) \subseteq
    \left( N_i \cup \{a_i^{(0)},\ldots,a_i^{(j)}\} \right) \setminus Y$ such that
	$\tricross_i^{(j)}(Y) \cap \overline{N_r} \neq \varnothing$, for all $r=1,\ldots,\ell$.
\end{subdividedLemma}
\begin{proof}
	(\ref{ThreeCrossExists})
    Observe that $S = \left( N_i \cup \{a_i^{(0)},\ldots,a_i^{(j)}\} \right) \setminus Y$ has cardinality
    $$|S|= k + (j + 1) - |Y| \ge k + 1.$$ 
    Applying property \ref{LargestCliqueSizeProperty} to $S$ produces a desired $3$-set 
    $T=\tricross_i^{(j)}(Y)$.
\end{proof}

A $3$-set $\tricross_i^{(j)}(Y)$, whose existence is established in 
Lemma \ref{ExistenceOfThreeCrosses}(\ref{ThreeCrossExists}), is 
called a {\em $3$-cross of $N_i^{(j)}$ with respect to $Y$}.
If $\tricross_i^{(j)}(Y)$ is a $3$-cross, then it is understood that 
$N_i \in {\cal N}^{(j)}$, $Y \subseteq N_i$, $|Y| \leq j$ and 
$\tricross_i^{(j)}(Y) \subseteq \left( N_i \cup \{a_i^{(0)},\ldots,a_i^{(j)}\} \right) \setminus Y$.

Now we enumerate several important properties of $3$-crosses.

\begin{subdividedLemma} \label{PrivateThreeCrossParts} 
	\item[] Suppose that $\tricross_i^{(j)}(Y)$ is 
    a $3$-cross.  Let $p = \tricross_i^{(j)}(Y) \cap N_i$.
    \item\label{ThreeCrossIntersectsHost} $|p| = 1$ or  $|p| = 2$.
    \item\label{PrivateElementProperty} If $|p| = 1$, then the vertex in $p$ is in $(N_i \setminus Y) \cap G$ and private to $N_i^{(j)}$. 
	\item\label{PrivatePairProperty} If $|p| = 2$, then $p$ is a $2$-set from $N_i \setminus Y$ and  
    $p$ is a private pair for $N_i^{(j)}$.
\end{subdividedLemma}

\begin{proof} 
	(\ref{ThreeCrossIntersectsHost})
	Because $|\tricross_i^{(j)}(Y)|=3$ and $\tricross_i^{(j)}(Y) \cap \overline{N_i} \neq \varnothing$,
	it follows that $|\tricross_i^{(j)}(Y) \cap N_i| \leq 2$.
	If  $|\tricross_i^{(j)}(Y) \cap N_i| = 0$, then $\tricross_i^{(j)}(Y) \subseteq  \{a_i^{(0)},\ldots,a_i^{(j)}\}$ which implies that,  because $\ell_j \geq 3$,
	there exists some $r$ such that $\tricross_i^{(j)}(Y) \subseteq N_r$, a contradiction.  Therefore, $1 \leq |p| \leq 2$.
	
	(\ref{PrivateElementProperty})  Because $|p| = 1$ and $\tricross_i^{(j)}(Y) \subseteq
    \left( N_i \cup \{a_i^{(0)},\ldots,a_i^{(j)}\} \right) \setminus Y$, it follows that
	$$\tricross_i^{(j)}(Y) = \{ a_i^{(\alpha)}, a_i^{(\beta)}, v\},$$
    for some $v \in N_i$ and $0 \leq \alpha < \beta \leq j$.  Note $p = \{v\}$ and $v \in N_i \setminus Y$.
	 
	Assume, to the contrary, that $v$ is not private to $N_i^{(j)}$. This means there exists 
	some $r \neq i$ such that $v \in N_r^{(j)}$.  Because $r\neq i$, the definition of the kernel implies
	$\{a_i^{(\alpha)}, a_i^{(\beta)}\} \subset N_r$;  therefore, $\tricross_i^{(j)}(Y) \subseteq N_r$.
	However, this implies $\tricross_i^{(j)}(Y) \cap \overline{N_r} = \varnothing$, 
	contradicting that $\tricross_i^{(j)}(Y)$ is a $3$-cross (Lemma \ref{ExistenceOfThreeCrosses}(\ref{ThreeCrossExists})).  
   So $v$ is private to $N_i^{(j)}$.  Because vertices in the kernel are double-covered (Lemma \ref{PairsInKernel}(\ref{PairsInKernelAreDoubleCovered})),
   it follows that $v \in G$.

	(\ref{PrivatePairProperty})
	We appropriately modify the argument given to establish (\ref{PrivateElementProperty}).  
	Because $|p| = 2$ and $\tricross_i^{(j)}(Y) \subseteq
    \left( N_i \cup \{a_i^{(0)},\ldots,a_i^{(j)}\} \right) \setminus Y$, it follows that
	$$\tricross_i^{(j)}(Y) = \{ a_i^{(\alpha)}, u, v\},$$
    for some $u, v \in N_i$ and $0 \leq \alpha \leq j$.  Note $p = \{u,v\}$ and $p \subseteq N_i \setminus Y$.
	 
	Assume, to the contrary, that $p$ is not a private pair for $N_i^{(j)}$. This means there exists 
	some $r \neq i$ such that $p \subset N_r^{(j)}$.  Because $r\neq i$, the definition of the kernel implies
	$a_i^{(\alpha)} \in N_r$;  therefore, $\tricross_i^{(j)}(Y) \subseteq N_r$.
	However, this implies $\tricross_i^{(j)}(Y) \cap \overline{N_r} = \varnothing$, 
	contradicting that $\tricross_i^{(j)}(Y)$ is a $3$-cross (Lemma \ref{ExistenceOfThreeCrosses}(\ref{ThreeCrossExists})).  
	So $p$ is private pair for $N_i^{(j)}$.
\end{proof}

\begin{subdividedLemma} \label{PrivatePairsForAll} 
	\item[] Suppose that $N_i \in {\cal N}^{(j)}$, $Y \subseteq N_i$ and $|Y| \leq j$.
    \item\label{PrivatePairWithNiceProperties} There exists a private pair, $p \subseteq N_i \setminus Y$, for $N_i^{(j)}$  
     such that either $p$ is a subset of a $3$-cross of $N_i^{(j)}$ with respect to $Y$, or $p$ does contain a private vertex for $N_i^{(j)}$.
\end{subdividedLemma}
\begin{proof} Lemma \ref{ExistenceOfThreeCrosses}(\ref{ThreeCrossExists}) guarantees there exists a
   $3$-cross $\tricross_i^{(j)}(Y)$.  Let $p = \tricross_i^{(j)}(Y) \cap N_i$; so $p \subset  N_i \setminus Y$.
   Lemma \ref{ExistenceOfThreeCrosses}(\ref{ThreeCrossIntersectsHost}) gives $|p| = 1$ or $|p| = 2$.
   
   If $|p| = 2$, then Lemma \ref{PrivateThreeCrossParts}(\ref{PrivatePairProperty}) shows
   that $p$ is the desired private pair for $N_i^{(j)}$.

   If $|p| = 1$, then Lemma \ref{PrivateThreeCrossParts}(\ref{PrivateElementProperty}) shows that
   $p$ contains a private vertex for $N_i^{(j)}$.  So it suffices now to note that
   a private pair for $N_i^{(j)}$ can be formed by adding any vertex from $N_i\setminus Y$ to $p$.
\end{proof}

Next we describe an inductive process to select a collection of private pairs for every $N_i$.
For each $i \in \{1,\ldots,\ell\}$, we define, by induction on time, a 
set $P_i = \left\{p_i^{(j)} : 0 \leq j \leq t_i \right\}$ of $t_i+1$ private pairs for $N_i$.  
The pair $p_i^{(j)}$ will be chosen from among the private pairs for $N_i$ existing at time $j$.  

Notice that, by property Lemma \ref{PairsInKernel}(\ref{PairsInKernelAreDoubleCovered}), 
any private pair contains at least one garbage vertex.  
Therefore, $p_i^{(j)}$ contains 
a vertex 
$g_i^{(j)} \in p_i^{(j)} \cap G$; call it the {\em anchor} of $p_i^{(j)}$.
The other vertex of $p_i^{(j)}$ is the {\em non-anchor}; it is denoted $u_i^{(j)}$.
It is possible that $\{g_i^{(j)},u_i^{(j)}\}\subset G$, but we shall still distinguish one vertex, $g_i^{(j)}$,
as the anchor of $p_i^{(j)}$.

We also define auxiliary sets $P_i^{(j)}=\left\{p_i^{(s)} : 0 \leq s \leq j\right\}$  and 
$G_i^{(j)}=\left\{g_i^{(s)} : 0 \leq s \leq j\right\}$. 
They are, respectively, the initial segments of the private pairs and the initial segments of the anchors for the private pairs selected for $N_i$ up to time $j$.  

Recall that every $k$-set $N_i$ has a final time, $t_i$, associated with it.  Now we associate a final time with every $g \in G$.
\begin{definition} \label{criticalTime} For $g \in G$, define the {\em critical time}
	for $g$, denoted $t_g$, to be the last stage at which there exists a $k$-set
	that contains $g$.  
	In other words,  $t_g = \max\left\{ t_i : g \in N_i \right\}.$
\end{definition}

We also need the following definition.

\begin{definition} \label{criticalIndices} For  $g \in G$, define the {\em critical index set}
	for $g$, denoted $I_g$, to be the set of indices of the $k$-sets that contain $g$ at stage $t_g$.  In other words,
	$$
	I_g=\{i :  g \in N_i, \; t_i=t_g\}.
	$$
\end{definition}

Observe that, by definition, $I_g \neq \varnothing$, for all $g \in G$. 
As we shall see in Lemma \ref{CriticalIndexSet}, the garbage vertices $g \in G$ with $|I_g|=1$ are particularly significant.  
The private pair selection process utilizes this information, so 
this motivates the following definition.

\begin{definition} \label{criticalGarbageVertices} For $x \in \{1,\ldots,\ell \}$, define {\em the set of critical garbage vertices}
	for $N_x$, denoted $\Gamma_x$, to be the set of $g \in N_x \cap G$ that have only $N_x$ in their critical index set.  In other words,
	$$
	\Gamma_x = \{ g \in G:  I_g = \{x\}\}.
	$$
\end{definition}
It is possible that $\Gamma_x = \varnothing$.  If $x_1 \neq x_2$, then by definition $\Gamma_{x_1} \cap \Gamma_{x_2} = \varnothing$.

Now we are ready to describe the private pair selection process.
Initially, for $i \in \{1,\ldots,\ell\}$, 
let $p_i^{(0)}= \left\{g_i^{(0)},u_i^{(0)}\right\} $ be a private pair for $N_i$ at time zero.
Such a private pair exists because of 
property Lemma \ref{PrivatePairsForAll}(\ref{PrivatePairWithNiceProperties}) with $Y=\varnothing$.
Set $G_i^{(0)} = \left\{ g_i^{(0)} \right\}$ and $P_i^{(0)} = \left\{ p_i^{(0)} \right\}$.

For $j>0$ and $i \in \{1,\ldots,\ell\}$, assume that the sets $P_i^{(j-1)}$ and $G_i^{(j-1)}$ have already been defined.  Also assume that a private pair $p_i^{(j)}=\left\{g_i^{(j)},u_i^{(j)}\right\}$ has already been chosen for each $N_i$ with $j \leq t_i$.
Now define 
\stack{$G_i^{(j)}$}{$=$}{$G_i^{(j-1)} \cup \left\{g_i^{(j)}\right\}$}{if $j \leq t_i$}{$G_i^{(j-1)}$}{if $j > t_i$,}
and similarly define,
\stack{$P_i^{(j)}$}{$=$}{$P_i^{(j-1)} \cup \left\{p_i^{(j)}\right\}$}{if $j \leq t_i$}{$P_i^{(j-1)}$}{if $j > t_i$.}
This definition yields $P_i^{(j)}=\left\{p_i^{(0)},\ldots,p_i^{(j)}\right\}$ and $G_i^{(j)}=\left\{g_i^{(0)},\ldots,g_i^{(j)}\right\}$;  
note that $\left|P_i^{(j)}\right| = \left|G_i^{(j)}\right| = j+1$, for each $0 \leq j \leq t_i$.

To complete the iterative process, it 
remains to describe how to select a private pair $p_i^{(j)}$, for each $i \in \{1,\ldots,\ell_j\}$. 

Firstly, and very importantly,  we prioritize the selection of a private pair for $N_i^{(j)}$ that contains a private vertex. Note that a private vertex of  any $k$-set belongs to $G$.  
If $N_i$ has a critical garbage vertex in $\Gamma_i \setminus G_i^{(j-1)}$ that is private to $N_i$ at time $j$, then
we select it. 
Otherwise, if $\Gamma_i \setminus G_i^{(j-1)}$ is empty,
but there is a private vertex in $N_i \setminus G_i^{(j-1)}$ that is private to $N_i$ at time $j$, we select it.
If either of these types of private vertex exists, then necessarily it is in $(N_i \cap G) \setminus G_i^{(j-1)}$; 
call it $g_i^{(j)}$.
Complete a private pair containing $g_i^{(j)}$ by adding a vertex
$a_x^{(j)}$, where $x \in \{1,\ldots,\ell_j\} \setminus\{i\}$.  Note that $a_x^{(j)} \in N_i$ by definition
of the kernel.
So, in this case, the pair selected is $p_i^{(j)} = \{g_i^{(j)}, a_x^{(j)}\}$, and it is a subset of $N_i$.
The fact that $p_i^{(j)}$ is a private pair for $N_i^{(j)}$ is due to $g_i^{(j)}$ being a private vertex for $N_i$ at time $j$.

If there is no vertex 
of $N_i \setminus G_i^{(j-1)}$ that is private to $N_i$ at time $j$,
apply Lemma \ref{PrivatePairsForAll}(\ref{PrivatePairWithNiceProperties}) 
with $Y=G_i^{(j-1)}$ to produce a pair  $\left\{g_i^{(j)},u\right\}$
single-covered by $\mathcal{N}^{(j)}$ such that  
$g_i^{(j)} \in G\setminus \left\{g_i^{(0)},\ldots,g_i^{(j-1)}\right\}$ and $u\in A\cup G$. 
If $u\in G$, then set $u_i^{(j)}=u$ which completes the private pair selection in this case. 
Otherwise, in anticipation of the need for well-behaved private pairs later, we  adjust 
$\left\{g_i^{(j)},g\right\}$. 
In this case $u \not\in G$, so we may assume that $u=a_r^{(s)}$ 
for some  $1\leq r\leq\ell$ and $0\leq s\leq t_r$.   Note that $r \neq i$ because
$a_r^{(s)} \in N_i$.   We now `promote' $a_r^{(s)}$, which means
setting the non-anchor vertex $u_i^{(j)}$ of the private pair to the latest possible kernel vertex that substitutes
for $a_r^{(s)}$ (preserving the private pair property).  This is accomplished as follows:

\begin{equation}
\label{prom}
u_i^{(j)}=\left\{
\begin{array}{cll}
a_r^{(s)}  & \hbox{if}\quad    j< s \leq t_r\\
\\
a_r^{(y)}  &  \hbox{if}\quad   s\leq  j\leq y, \hbox{\; where } y=\min\{t_i,t_r\}\\
\\
a_x^{(t_i)}  & \hbox{if}\quad    s\leq t_r<j\leq t_i,  \hbox{\; with any}\; x\in\{1,\ldots,\ell_{t_i}\}\setminus\{i\} .
\end{array}
\right.
\end{equation}

Notice that promotion is not applied if $g_i^{(j)}$ is a  private vertex. Furthermore, (\ref{prom}) changes the non-anchor (that is $u_i^{(j)} \neq u$) only if $s \leq j$.
Also observe that if $u_i^{(j)} = a_r^{(y)}$ with $y \not\in \{t_i,t_r\}$, 
then the 
non-anchor is unchanged (that is, $u_i^{(j)} = u$ under these conditions).
We next prove that this promotion process in fact produces a private pair  for $N_i$ at time $j$, furthermore, anchor vertices are not repeated for the private  pairs in $N_i$.

\begin{lemma}
\label{ppp}  
The pair $p_i^{(j)}= \left\{g_i^{(j)},u_i^{(j)}\right\}$ as defined by {\em (\ref{prom})} 
 is a private pair for $N_i$ at time $j$. 
Furthermore, $g_i^{(j)} \in G\setminus \left\{g_i^{(0)},\ldots,g_i^{(j-1)}\right\}$, $j=1,\ldots,t_i$, and if
$u_i^{(j)}$ is a kernel vertex, then 
$u_i^{(j)}=a_x^{(y)}$, where  $j\leq y$ and $x\leq \ell$.
\end{lemma}
\begin{proof}  We may assume that
no vertex 
of $N_i \setminus G_i^{(j-1)}$ is private to $N_i$ at time $j$, since otherwise, the private pair selected has the form
$p_i^{(j)} = \{g_i^{(j)}, a_x^{(j)}\}$
  satisfying the claim with $j=y$.
 
Notice that the anchor of $p_i^{(j)}$ is not affected by the promotion process,
so $g_i^{(j)} \in G\setminus \left\{g_i^{(0)},\ldots,g_i^{(j-1)}\right\}$
follows from the generation of $\left\{g_i^{(j)},u\right\}$ via 
Lemma \ref{PrivatePairsForAll}(\ref{PrivatePairWithNiceProperties}). 
Assuming that the non-anchor is a kernel vertex, the formula (\ref{prom}) yields either 
$u_i^{(j)}=a_r^{(y)}$, where $y=s$ and $j<s$, or $u_i^{(j)}=a_x^{(y)}$, where $y\in\{t_i,t_r\}$; 
in each case  $j\leq y$.    
Thus, the second part of the claim follows.

It remains to show that 
$p_i^{(j)}= \left\{g_i^{(j)},u_i^{(j)}\right\}$, as defined by (\ref{prom}), is a private pair for $N_i^{(j)}$.
First we verify that $p_i^{(j)} \subset N_i$.  Clearly, $g_i^{(j)} \in N_i$, so it suffices to show that promotion
produces $u _i^{(j)} \in N_i$.
If promotion does not change the non-anchor (line 1 of (\ref{prom})), then
$a_r^{(s)} \in N_i$ follows from the generation of $\left\{g_i^{(j)},u\right\}$ via 
Lemma \ref{PrivatePairsForAll}(\ref{PrivatePairWithNiceProperties}). 
So we may assume that promotion does change $u$.
Because $i \neq r$, if $t_r < t_i$, then $a_r^{(t_r)} \in N_i$; whereas
if $t_i \leq t_r$, then $a_r^{(t_i)} \in N_i$ (line 2 of (\ref{prom})).  Similar reasoning shows
that $x \neq i$ implies $a_x^{(t_i)} \in N_i$ (line 3 of (\ref{prom})) because $x\in\{1,\ldots,\ell_{t_i}\}\setminus\{i\}$ 
means $t_i \leq t_x$.

Now we verify the privacy of $p_i^{(j)}$.
If $u_i^{(j)}$ is set to $a_r^{(s)}$ (that is, $j < s$  corresponding to line 1 of (\ref{prom})), then promotion does not change the non-anchor which means $p_i^{(j)}$
retains the privacy granted from the application of Lemma \ref{PrivatePairsForAll}(\ref{PrivatePairWithNiceProperties}) which 
generated it.  
So we may assume that promotion does change the non-anchor; that is, $s \leq j$.  
Because $\{g_i^{(j)},a_r^{(s)}\}$ is single-covered by $N_i$ at time $j$ and $s \leq j$, the anchor vertex 
$g_i^{(j)}$ does not belong to any
 $N_x\in \mathcal{N}^{(s)}\setminus\{N_r\}$ such that $t_x\geq j$. In particular, 
 if $\{g_i^{(j)},a_x^{(t_i)}\}\subset N_q$, and $N_q\in\mathcal{N}^{(j)}$, then $q\in\{i,r\}$. 
Hence the pair $\{g_i^{(j)},a_x^{(t_i)}\}$ is single-covered
by $N_i$ at time $j$, if either $x=r$ or $t_r<j$. These cases correspond  to the second and third line of (1).

In the exceptional case when $j\leq t_r<t_i$  the second line of (1) yields $p_i^{(j)}=\{g_i^{(j)},a_r^{(t_r)}\}$. This is a single-covered pair, since $s\leq t_r$ and $\{g_i^{(j)},a_r^{(s)},\}$ is single-covered  at time $j$.
\end{proof}

The next lemma establishes a property of the selected private pairs
that is essential to building small transversal sets for arcs of the digraph $D$ constructed in Section \ref{digraph}.

\begin{lemma} \label{PPinsideThreeCross}
Suppose that $p_i^{(j)}= \left\{ g_i^{(j)},u_i^{(j)}\right\}$ is a selected private pair for $N_i^{(j)}$.
If $u_i^{(j)}=a_x^{(y)}$ where $y \not\in \{t_i, t_x\}$ and $g_i^{(j)}$ is not a private vertex for $N_i$ at time $j$,
then $p_i^{(j)}$ is a subset of a $3$-cross of $N_i^{(j)}$ with respect to $G_i^{(j-1)}$.
\end{lemma}
\begin{proof}  Because the private pair selection process prioritizes
   selecting private pairs with private vertices, the hypothesis that $g_i^{(j)}$ is not a private vertex for $N_i$ at time $j$
   means $N_i \setminus G_i^{(j-1)}$ has no private vertices at time $j$.
   Consequently, the process to select $p_i^{(j)}$ must first have applied
   Lemma \ref{PrivatePairsForAll}(\ref{PrivatePairWithNiceProperties}) to generate a pair
   $\left\{g_i^{(j)},u\right\}$ in which $u$ is a kernel vertex (since the final private pair has a non-anchor in the kernel).  
   Lemma \ref{PrivatePairsForAll}(\ref{PrivatePairWithNiceProperties}) guarantees
   a $3$-cross $\tricross_i^{(j)}(G_i^{(j-1)})$ containing $\left\{g_i^{(j)},u\right\}$.   The promotion process must not 
   have changed the non-anchor $u$ because the final non-anchor is $u_i^{(j)}= a_x^{(y)}$ with $y \not\in \{t_i, t_x\}$ (see observation prior to Lemma
   \ref{ppp}).
   Therefore, $p_i^{(j)}= \left\{g_i^{(j)},u\right\} \subseteq \tricross_i^{(j)}(G_i^{(j-1)})$, as desired.
\end{proof}

Let $P^{(j)} = \left\{p_i^{(j)} : 1 \leq i \leq \ell_j \right\}$ denote the private pairs defined by this process at time $j$ and let $P_i = \left\{p_i^{(j)} : 0 \leq j \leq t_i \right\}$ be the set of $t_i+1$ private pairs defined for $N_i$ by this process.  The collection of all 
selected pairs is defined as
$$
P = \bigcup_{j=0}^t P^{(j)} =\bigcup_{i=0}^\ell P_i.
$$ 
\noindent

Next we list for reference obvious properties of the private pairs in $P$ which are summarized in Lemma 8 of \cite{SzemerediPetruska}:
\begin{itemize}
\item[-] $P^{(j_1)} \cap P^{(j_2)} = \varnothing$, for $0 \leq j_1 < j_2 \leq t$;
\item[-] any pair in $\bigcup_{s=0}^j P^{(s)}$ is at most single-covered by ${\cal N}^{(j)}$, for $j = 0,\ldots,t$;
\item[-] $\left|P^{(j)}\right| = \ell_j$ and every $\left|P_i \cap P^{(j)}\right| = 1$, for all $0 \leq j \leq t$ and $1 \leq i \leq \ell_j$. 
\end{itemize}

\section{Advancement}
\label{advancement}

In this section we describe a technical modification to the private pairs, called {\em advancement}, that enables our final tight asymptotic bounds.
This modification concerns rare, but troublesome private pairs that we alter and reorder.  During this reordering, troublesome private pairs 
are given improved non-anchors and advanced in time (hence the process name) to protect their anchors, guaranteeing that the anchors are swapped
via the  recursive process (\ref{RecursiveDefinitionOfMs})  defined later in Section \ref{Skew}.
Advancement preserves the number of pairs in $P$, but a few pairs will lose their private nature, sacrificing themselves to advance vital
troublesome private pairs.
To limit the size of the required sacrifices and identify these vital troublesome pairs, we first 
partition the $k$-sets into three types :
a $k$-set $N_x$ is {\em weightless} if $|\Gamma_x| = 0$; it is {\em light}
if $0 < |\Gamma_x| < m^{1/3}$ and {\em heavy} if $ m^{1/3} \leq |\Gamma_x|$.  

The justification for this otherwise arbitrary appearing trichotomy is that 
it optimizes the bounds in Theorem \ref{maintheorem}.  As we shall see, weightless $k$-sets are innocuous, so there will be no need
to address them further.  Heavy sets require too many sacrifices to correct via advancement; their challenges
will be settled efficiently by observing that there are few heavy $k$-sets (see Lemma \ref{boundHeavyKSets})
so a crude solution is affordable (definition of $T_H$ in Section \ref{transversal}).  This leaves 
the correction of light $k$-sets that define the `troublesome' pairs (see Definition \ref{troublesome}).  

Before defining troublesome pairs, we define and examine this precursor to them:

\begin{definition} \label{problematic}
  A private pair $p_i^{(j)} \in P$ is {\em problematic} if it has the form $p_i^{(j)} = \{g,a_x^{(y)} \}$ with $j \leq y < t_x$ and $g \in N_x $. 
\end{definition}

\noindent
The following lemmas present important properties of problematic private pairs.

\begin{lemma}
	\label{CriticalIndexSet}
	If $p_i^{(j)} = \{g,a_x^{(y)}\}$ is a problematic private pair, then
	\begin{enumerate}[label=(\roman*)]
		\item \label{HeadTailKSetsDistinct} $i \neq x$,
		\item \label{SparseInterval} for any $z \in \{1,\ldots,\ell\}$, if $g \in N_z$ and $y \leq t_z$, then $z \in \{i,x\}$,
		\item \label{CriticalIndexSetsOfArcs} $I_g \subseteq \{i,x\}$,
	\end{enumerate}
\end{lemma}
\begin{proof}
	\ref{HeadTailKSetsDistinct}
	Because $a_x^{(y)} \in p_i^{(j)} \subseteq N_i$ and $a_x^{(y)} \not\in N_x$, it follows that $x \neq i$.
	
	\ref{SparseInterval}
	Assume, to the contrary, that 
	there exists $z$ such that $g \in N_z$, $t_z \geq y$, and $z \not\in \{i,x\}$.  
	Because $y \leq t_z$ and $z \neq x$, it follows that $a_x^{(y)} \in N_z$. 
	Because $j \leq y$, we have  
	$p_i^{(j)} = \{g,a_x^{(y)}\} \subseteq N_i^{(j)} \cap N_z^{(j)}$,
	contradicting that $p_i^{(j)}$ is private for $N_i$ at stage $j$. 
	
	\ref{CriticalIndexSetsOfArcs}
	Consider $z \in I_g$; so $g \in N_z$ and $t_z=t_g$.  Because $g \in N_x$, we have $t_x \leq t_g$. Hence, $y < t_x \leq t_z = t_g$.  
	Now \ref{SparseInterval} implies that $z \in \{i,x\}$.  Therefore, $I_g \subseteq \{i,x\}$.
\end{proof}

Lemma \ref{CriticalIndexSet}\ref{HeadTailKSetsDistinct} shows  that the anchor of a problematic private pair is not a private vertex, thus  the pair is obtained through promotion (\ref{prom}). The next lemma highlights the size of a critical index set of the anchor vertex,  in particular, distinguishing whether
it is one or two.
\begin{lemma} \label{Trichotomy}
	Suppose $p_i^{(j)} = \{g,a_x^{(y)}\}$ is a problematic private pair.
	\begin{enumerate}[label=(\roman*)]
		\item \label{Tri} $|I_g| = 1$ or $|I_g| = 2$.
		\item \label{TypeW} If $|I_g| = 2$, then $I_g = \{i,x\}$ and there are at most two problematic private pairs with anchor $g$.
		\item\label{DifficultCase} For $|I_g| = 1$, if $I_g = \{i\}$ then $p_i^{(j)}$ is the only problematic private pair with anchor $g$; otherwise, $I_g = \{x\}$ and every problematic private pair with anchor $g$ has a non-anchor from $\overline{N_x}$.
	\end{enumerate}
\end{lemma}
\begin{proof}  \ref{Tri} Because $I_g \neq \varnothing$, this claim is just a rephrasing of Lemma \ref{CriticalIndexSet}\ref{CriticalIndexSetsOfArcs}.

	\ref{TypeW}
	Because $|I_g| = 2$ and Lemma \ref{CriticalIndexSet}\ref{CriticalIndexSetsOfArcs}, it follows that
	$I_g = \{i,x\}$.  	Consider an arbitrary problematic private pair $p_u^{(\alpha)}$ with anchor $g$.  Lemma \ref{CriticalIndexSet}\ref{CriticalIndexSetsOfArcs}
	implies $u \in \{i,x\}$.
	Now simply observe that $N_i$ and $N_x$ each may have at most one private pair containing $g$ because anchors are never repeated in the private pair selection process. 		
	
	\ref{DifficultCase} 
	By Lemma \ref{CriticalIndexSet}\ref{CriticalIndexSetsOfArcs}, $\{z\} = I_g \subset \{i,x\}$.
	If $z = i$, then $p_i^{(j)}$ is the unique private pair for $N_i$ with anchor $g$;
	otherwise, $z = x$.  Consequently, every problematic private pair with anchor $g$ that is not a private pair for $N_x$
	has non-anchor from $\overline{N_x}$.
\end{proof}

Now we are ready to define troublesome private pairs.  

\begin{definition} \label{troublesome}
	A private pair $p_i^{(j)} = \{g, a_x^{(y)}\}$ is {\em troublesome} if $p_i^{(j)}$ is problematic, $N_i$ is light, and $g \in \Gamma_i$.
\end{definition}

These pairs are a concern because they cause particularly thorny conflicts that must be removed to 
obtain a skew $(2,m)$-system that appears in Theorem \ref{SkewSystem}.
The next lemma introduces a possible replacement pair for a troublesome pair.  The {\em replacement pair} for $p_i^{(j)}$ is a new pair $\hat{p}_i^{(j)}$ with the same anchor but
improved non-anchor, $a_x^{(t_x)}$, that is immune to swaps (see the recursive process (\ref{RecursiveDefinitionOfMs}) in Section \ref{Skew}) because $t_x$ is a terminal stage for $N_x$.  
This immunity ensures that, in the ultimate skew $(2,m)$-system considered, the replacement pair will continue to intersect $m$-sets that contain $a_x^{(t_x)}$.
Note however that, in contrast to the original pair, the replacement pair may not be private for $N_i$ at time $j$, though it is private for $N_i$ at later stages, as the next  lemma specifies.

\begin{lemma} \label{replacementPair}
	If $p_i^{(j)} = \{g, a_x^{(y)}\}$ is a troublesome private pair, then the replacement pair $\hat{p}_i^{(j)} = \{g, a_x^{(t_x)}\}$ is
	private for $N_i^{(\alpha)}$, for any $\alpha \geq y$.
\end{lemma}
\begin{proof}
Lemma \ref{CriticalIndexSet}\ref{HeadTailKSetsDistinct} says $x \neq i$,  and $t_x < t_i$ because $g \in N_x$ and $I_{g} = \{i\}$; therefore, $a_x^{(t_x)} \in N_i$.  This means that the pair $\{g, a_x^{(t_x)}\}$ is a subset of $N_i$.
Indeed, Lemma \ref{CriticalIndexSet}\ref{SparseInterval} implies that this pair is private for $N_i^{(y)}$.
Consequently, this pair is private for $N_i$ at any time $\alpha$ satisfying $\alpha \geq y$.
\end{proof}

Because a replacement private pair is not necessarily private at the same time as the original, to use this new pair we
must advance the replacement pair in time at the expense of destroying another private pair;
this is ``advancement''.  The advancement expense is justified because these replacement pairs contain important anchors 
whose swap must be protected to prevent many other conflicts.   
Only light $k$-sets are considered so that the expense of advancing these replacement pairs is limited.   
In some cases it may not be possible to advance a given troublesome pair because no pair remains to be sacrificed.
The process takes this possibility into account (see `neutralized' pairs below).

\vskip 0.15cm
\noindent
{\sc Advancement}:
\noindent
Here is the formal description of the advancement process.
For each light $k$-set $N_i$, 
we view the collection $P_i$ of private pairs as a list:
$p_i^{(0)},\ldots, p_i^{(t_i)}$.
The process treats troublesome private pairs in this list according to the order in which they appear {\textemdash} earlier troublesome pairs are processed before later ones.
Suppose that $p_i^{(j)} = \{g_i^{(j)}, a_x^{(y)}\}$ is the next troublesome pair to consider for $N_i$.
By definition, this means $g_i^{(j)} \in \Gamma_i$, that is $g_i^{(j)}$ is a critical garbage  vertex in $N_i$.  
Because  $p_i^{(j)}$ is problematic, $j\leq y < t_x$.

If there exists a private pair $p_i^{(\alpha)} =  \{g_i^{(\alpha)}, u_i^{(\alpha)}\}$ for $N_i$ with $\alpha \geq y$
and $g_i^{(\alpha)} \not\in \Gamma_i$, then set $\{g_i^{(j)}, a_x^{(t_x)}\}$ as the ``replacement pair'' to be the new private pair for $N_i$ at time $\alpha$, and 
move the pair $\{g_i^{(\alpha)}, u_i^{(\alpha)}\}$ to time $j$ to take the place of the displaced original troublesome pair.  
In this case we refer to the replacement pair $\{g_i^{(j)}, a_x^{(t_x)}\}$, which is the new $p_i^{(\alpha)}$, as {\em advanced}; 
the moved pair  $\{g_i^{(\alpha)}, u_i^{(\alpha)}\}$ is labeled {\em deactivated} as it may no longer be private to $N_i$ at the stage at which it now appears, time $j$.
If there
is no private pair $p_i^{(\alpha)} =  \{g_i^{(\alpha)}, u_i^{(\alpha)}\}$ for $N_i$ with $\alpha \geq y$
and $g_i^{(\alpha)} \not\in \Gamma_i$, then we simply label the troublesome pair $p_i^{(j)} = \{g_i^{(j)}, a_x^{(y)}\}$ {\em neutralized}.

After all troublesome pairs for $N_i$ are processed, some pairs are deactivated, and the remaining pairs are all private pairs
for $N_i$ at the times at which they appear since Lemma \ref{replacementPair} guarantees that replacement pairs are private.
The only remaining troublesome pairs are neutralized.  Also note that this advancement process does not alter the set of anchors.  
In particular, distinct private pairs for $N_i$ still have different anchors.  

To clearly indicate post-advancement definitions, we employ math bold face font for these new sets.  For example,
the new set of pairs for $N_i$ is denoted $\mathbb{P}_i$.  The union of all new pairs for all $N_i$'s 
is denoted $\mathbb{P}$.
Similarly, 
$\mathbb{G}_i^{(j)}$ are anchors for the pairs in $\mathbb{P}_i$ for $N_i$ up to time $j$. 
A pair in $\mathbb{P}$ that is not deactivated is {\em active}.  Observe that active pairs are private pairs.

\section{A digraph} 
\label{digraph}

 In  this section we introduce an auxiliary digraph $D$ on the kernel $A$ as its vertex set, where
 a vertex $a_x^{(y)}$  represents the private pair $p_x^{(y)} \in \mathbb{P}$. The arcs of $D$ will be labeled with garbage vertices and will serve as the code of `conflicting' private pairs. 
Thus an independent set in $D$ will determine a special subset of non-conflicting private pairs 
 in $\mathbb{P}$.  These pairs represented by the vertices of a large enough independent set of $D$ will be used in Section \ref{Skew} to define a large skew cross-intersecting $(2,m)$-system.

The final   skew cross-intersecting $(2,m)$-system will be generated by the recursive process (\ref{RecursiveDefinitionOfMs})  in Section \ref{Skew}. To protect critical garbage vertices during this process we introduce the following definition.

\begin{definition} \label{protectedPrivateAnchors} For $x \in \{1,\ldots,\ell \}$, define the set of {\em protected private anchors for $N_x$ at time $y$}, 
	denoted $\Lambda_x^{(y)}$, to be the set:
	$$
	\Lambda_x^{(y)} = \left( \mathbb{G}_x^{(y-1)} \cap \Gamma_x \right) \setminus \{ g \in G: 
	\mbox{ $g$ is the anchor of a  neutralized private pair of $N_x$}\}.
	$$
\end{definition}
This set, $\Lambda_x^{(y)}$, contains all anchors for non-neutralized private pairs for $N_x$ up to time $y$.
In the special case in which $N_x$ is light, these anchors will
get swapped via the recursive procedure (\ref{RecursiveDefinitionOfMs}) to produce $M_x^{(y)}$ (defined later in Section \ref{Skew}); this explains
condition \ref{LightKSetException} below.

With these definitions complete, we are now ready to define an arc-labeled digraph $D$ on the vertex set $A$.
Recall that $A$ is the set of kernel vertices; and  every element of $A$ has the form
$a_x^{(y)}$, where $1 \leq x \leq \ell$ and $0 \leq y \leq t_x$.  
The arcs of the digraph $D$ are defined as follows.
\begin{equation} \label{arcDefinition} 
   a_i^{(j)} \stackrel{g}{\rightarrow} a_x^{(y)} \mbox{ is an arc of } D 
   \end{equation}
   if and only if all the following are true:
   \begin{quote}
   \begin{enumerate}[label=(\ref{arcDefinition}.\arabic{enumi})] 
       \item \label{PrivatePairFromArc} $p_{i}^{(j)} = \{ g, a_x^{(y)}\} \in \mathbb{P}_i$ and $p_{i}^{(j)}$ is active, 
       \item \label{AnchorArcProperty} $g \in N_x \cap N_i$,  
       \item \label{IndexArcProperty} $j \leq y < t_x$, and  if $j = y$ then $i < x$,
       \item \label{LightKSetException} if $N_x$ is light, then $g \not\in \Lambda_x^{(y)}$.
   \end{enumerate}
   \end{quote}
\noindent

Together, conditions \ref{PrivatePairFromArc}, \ref{AnchorArcProperty}, \ref{IndexArcProperty} imply
that $p_{i}^{(j)}$ is a problematic pair.  Because conditions \ref{PrivatePairFromArc} requires $p_{i}^{(j)}$ is active,
the pair is also a private pair for $N_i$.  Condition \ref{IndexArcProperty} specifies that $y < t_x$ which means
that $p_{i}^{(j)}$ is not an advanced pair.

We refer to the vertex $a_i^{(j)}$ as the {\em tail} of the arc $a_i^{(j)} \stackrel{g}{\rightarrow} a_x^{(y)}$;
naturally,  $a_x^{(y)}$ is the {\em head}.
Observe that \ref{PrivatePairFromArc} implies $g \in G$ because the private pair $p_{i}^{(j)}$ must intersect $G$.
It also implies the digraph $D$ has out-degree at most one because an out-going arc from vertex $a_i^{(j)}$, if there is one,
is determined by $p_{i}^{(j)}$.
Condition \ref{AnchorArcProperty} guarantees that $g$ is not a private vertex to $N_i$ at time $j$;
that is, the production of the private pair $p_{i}^{(j)} = \{ g, a_x^{(y)}\}$ must have invoked 
Lemma \ref{PrivatePairsForAll}(\ref{PrivatePairWithNiceProperties}).

Lemma \ref{Trichotomy}
shows that if there are more than two arcs in $D$ with label $g$, then $|I_g| = 1$.  
This explains the focus on such vertices given in Definition \ref{criticalGarbageVertices}.
  
The idea is that a labeled arc in $D$ is a code for 
conflicting elements of an initial $(2,m)$-system developed in Section \ref{Skew}.  
By eliminating the conflicts
exposed by these arcs, this $(2,m)$-system will eventually be reduced to a valid skew cross-intersecting $(2,m)$-system (Theorem \ref{SkewSystem}). 
The first three arc-defining conditions
naturally encode the conflicts, but
condition \ref{LightKSetException} is artificial and technical, arising from the special manner in which we
treat light $k$-sets in Section \ref{transversal}.
Conflict elimination is achieved via a large independent set in $D$, so we now turn 
our attention to guaranteeing such a set.

Given an independent set $F$ of $D$, call a pair  $p_i^{(j)} \in \mathbb{P}$ {\em free} if $a_i^{(j)} \in F$. 
The size of $F$ determines this size of the final skew cross-intersecting $(2,m)$-system.
The complement $T=A \setminus F$ is a transversal (a vertex cover) of the arcs in $D$.
We strive for an upper bound on $T$, which yields a lower bound on $F$ and thus on the number of free pairs.
In the next section we shall determine a small transversal.
The selection of this transversal will apply the aforementioned trichotomy of $k$-sets, a concept 
originating from the Lemma \ref{Trichotomy}. 

\section{A small transversal set} 
\label{transversal}

Our goal in this section  is to find a large independent set in digraph $D$ defined in Section \ref{digraph}, or equivalently, a small transversal set (vertex cover) of the arcs of $D$.
More specifically, we seek a transversal set $T \subset A$ such that $|T| = O(m^{5/2})$.
This will guarantee a set of free pairs, $F=A \setminus T$, that has cardinality $|F| \geq |A| - O(m^{5/2})$
which will play an essential role in the proof of Theorem \ref{maintheorem}.
We build $T$ as the union of four subsets $T_W$, $T_H$, $T_S$, and $T_L$ that we describe below.  Each of these
subsets targets arcs of $D$ that arise from different circumstances.

Since the tail of an arc from $D$ represents a problematic private pair,  Lemma \ref{CriticalIndexSet}
shows that if an arc in $D$ has the label $g$, then $|I_g| = 2$ or $|I_g| = 1$.
All arcs arising from the case $|I_g| = 2$  will be covered by the set $T_W$ of their tails:
$$T_W=\{a_i^{(j)}: a_i^{(j)} \stackrel{g}{\longrightarrow} a_x^{(y)} \hbox{ is an arc of $D$ and $|I_g| = 2$}\}.$$
\noindent
The next lemma proves that $T_W$ is not too large.
\begin{lemma}  \label{SmallWTransversal}
	$|T_W| \leq 6m$
\end{lemma}
\begin{proof}
	Lemma \ref{CriticalIndexSet} shows that if  $|I_g| = 2$, then there are at most two arcs in $D$ labeled $g$.
Because	$|G| \leq 3m$ (Lemma \ref{FirstLemma}(\ref{GUpperBound})), it follows that there are at most $2\cdot 3m$ arcs of this type.
\end{proof}

We next turn our attention to covering arcs in $D$ that have a label $g$ such that $|I_g| = 1$.
Such arcs have a head vertex in a light or heavy $k$-set.  We address the latter kind first.
We form a small transversal for these arcs using their head vertices (see definition of $T_H$ below).  First we prove
that there are not many heavy $k$-sets.

\begin{lemma} \label{boundHeavyKSets}
	The number of heavy $k$-sets is at most $3m^{2/3}$.
\end{lemma} 
\begin{proof}  Let $h$ denote the number of heavy $k$-sets.
	By definition, a $k$-set $N_x$ is heavy if $ m^{1/3} \leq |\Gamma_x|$.  Because
	different heavy sets have disjoint sets of critical garbage vertices, we have
	$$
	|G| \geq \biggl| \bigcup_{\tiny N_x \mbox{ heavy }} \Gamma_x \biggl|
	= \sum_{\tiny N_x \mbox{ heavy }} | \Gamma_x | \geq h \cdot m^{1/3}
	$$
Now $h \leq  3m^{2/3}$ follows from $|G| \leq 3m$ that is a consequence of Lemma \ref{FirstLemma}(\ref{GUpperBound}).
\end{proof}

Now define the transversal vertices contributed by the heavy $k$-sets:
$$T_H=\{a_x^{(y)} : a_i^{(j)} \stackrel{g}{\longrightarrow} a_x^{(y)} \hbox{ is an arc of $D$ and $N_x$ is heavy}\}.$$
The next lemma proves that this set is not too large.

\begin{lemma} \label{SmallHTransversal}
	$|T_H| \leq 3m^{5/3}$
\end{lemma}
\begin{proof}
	A heavy $k$-set $N_x$ contributes at most $t_x \leq t$ elements to $T_H$.
	Because $t < m$ (Lemma \ref{FirstLemma}(\ref{tUpperBound})) and the number of heavy $k$-sets is at most $3m^{2/3}$ (Lemma \ref{boundHeavyKSets}),
	we conclude that $|T_H| \leq m \cdot 3m^{2/3}.$
\end{proof}

Finally, we address the remaining arcs in $D$; these have a label $g$ such that $|I_g| = 1$ and their head vertex is in a light $k$-set.
Completing  the final construction of our small transversal set $T$ is a very important property we seek to ensure, Property~\eqref{property:L}:
\begin{equation}
	\tag{L}\label{property:L}
	\parbox{\dimexpr\linewidth-4em}{%
		\mathstrut
		If $N_x$ is a light $k$-set, $g$ is the anchor for $p_x^{(y)}$, and $g\in \Gamma_x$, \\
		then either $p_x^{(y)}$ is neutralized or $a_x^{(y)} \not\in T$.
		\mathstrut
	}
\end{equation}
\noindent  
The conclusion	$a_x^{(y)} \not\in T$ is equivalent to $a_x^{(y)} \in F$ or $p_x^{(y)}$ is free.
Property~\eqref{property:L} guarantees that, if $p_x^{(y)}$ is not neutralized, then $g$ gets swapped into $M_x^{(y+1)}$ in the recursive process (\ref{RecursiveDefinitionOfMs})
generating the final large skew cross-intersecting $(2,m)$-system in Section \ref{Skew}, thereby avoiding many conflicts present in the initial
system.
Guaranteeing Property~\eqref{property:L} for heavy $k$-sets seems to require a large transversal set, frustrating the primary 
small transversal objective. 
This explains the name of Property~\eqref{property:L}; it only concerns light $k$-sets.

Observe that $T_H$ contains the head vertices of all arcs in $D$ that fall into a heavy $k$-set.
In contrast, we use tail vertices to cover arcs in $D$ whose label $g$ satisfies $|I_g| = 1$ and whose head falls into a light $k$-set.  
We do this to satisfy Property~\eqref{property:L}.  In particular,
define these two sets
\begin{eqnarray*}
	T_S=\{a_i^{(j)} \in A \setminus T_W: & \hbox{$p_i^{(j)}$ is deactivated, or } \hfill \\ 
	\ & a_i^{(j)} \stackrel{g}{\longrightarrow} a_x^{(y)} \hbox{ is an arc of $D$ and $p_i^{(j)}$ is neutralized}\},
\end{eqnarray*}

and
$$T_L=\{a_i^{(j)} \in A \setminus (T_W \cup T_S): a_i^{(j)} \stackrel{g}{\longrightarrow} a_x^{(y)} \hbox{ is an arc of $D$ and $N_x$ is light}\}.$$
\noindent
Observe that, in both definitions,
the condition $a_i^{(j)} \in A \setminus T_W$ guarantees that label $g$ satisfies $|I_g| = 1$.

The final choice of transversal set is: 
$$
T = T_W \cup T_H \cup T_S \cup T_L.
$$
Clearly $T$ is a transversal for the arcs of $D$.  We now turn to proving that $T$ satisfies Property~\eqref{property:L}
and $|T| = O(m^{5/3})$.  We consider Property~\eqref{property:L} first. 

\begin{lemma} \label{PropertyLSatisfied}
	$T$ satisfies Property~\eqref{property:L}
\end{lemma}
\begin{proof}
	Suppose $N_i$ is a light $k$-set, $g$ is the anchor for $p_i^{(j)}$, and $g \in \Gamma_i$.  Assume that $p_i^{(j)}$ is not neutralized. We must prove
	that $a_i^{(j)} \not\in T$.  First observe that $a_i^{(j)} \not\in T_W \cup T_H$ by definition of $T_W$ and $T_H$, because $|I_g|\neq 2$ and $N_i$ is a light $k$-set. 
	
	Next we claim there is no arc of $D$ in which $a_i^{(j)}$ is the tail.  Suppose, to the contrary, that
	$a_i^{(j)} \stackrel{g}{\rightarrow} a_x^{(y)}$ is such an arc of $D$.  So, by condition \ref{PrivatePairFromArc}, we have $p_i^{(j)}=\{g,a_x^{(y)}\}$.
	Because $p_i^{(j)} \subseteq N_i$, clearly $x \neq i$ since otherwise $a_i^{(y)} = a_x^{(y)} \in N_i$, a contradiction.
	Recall that, if	$a_i^{(j)} \stackrel{g}{\rightarrow} a_x^{(y)}$ is an arc of $D$, then  $p_i^{(j)}$ is a problematic pair.
	Since we additionally assume that $g \in \Gamma_i$ and $N_i$ is light, it follows that $p_i^{(j)}$ is a troublesome pair.  However, we have also assumed that $p_i^{(j)}$ is not neutralized.
	Therefore, $p_i^{(j)}$ is an advanced troublesome pair; that is, $p_i^{(j)}=\{g,a_x^{(t_x)}\}$.  In particular, $y = t_x$ which contradicts condition \ref{IndexArcProperty} that
    $a_i^{(j)} \stackrel{g}{\rightarrow} a_x^{(y)}$ is an arc of $D$.
	
	Because  $p_i^{(j)}$ is not neutralized, $a_i^{(j)} \not\in T_S$.  We noted earlier that $a_i^{(j)} \not\in T_W \cup T_H$.   The last paragraph shows that $p_i^{(j)}$ not neutralized implies there is no arc of $D$ that has $a_i^{(j)}$ as a tail.
	It follows that $a_i^{(j)} \not\in T_L$ since these sets are defined taking only tails of arcs.   We conclude that $a_i^{(j)} \not\in T$, as claimed.
\end{proof}

To prove $|T| = O(m^{5/3})$, Lemma \ref{SmallWTransversal} and Lemma \ref{SmallHTransversal} imply
that it suffices to prove that $|T_S \cup T_L| = O(m^{5/3})$.  
First we bound the number of light $k$-sets.

\begin{lemma} \label{boundLightKSets}
	The number of light $k$-sets is at most $3m$.
\end{lemma}
\begin{proof}
	Let $N_{\beta_1},\ldots,N_{\beta_q}$ be an enumeration of all the light $k$-sets in ${\cal N}$.
	Because $\varnothing \neq \Gamma_{\beta_i} \subseteq G$ for all $1 \leq i \leq q$ and these sets of critical garbage vertices
	are disjoint, we conclude that 
	$$
	q \leq \sum_{i=1}^{q} |\Gamma_{\beta_i}| = \bigg| \bigcup_{i=1}^{q} \Gamma_{\beta_i} \bigg| \leq |G|.
	$$ 
	Therefore $q \leq 3m$ because Lemma \ref{FirstLemma}(\ref{GUpperBound}) implies $|G| \leq 3m$.
\end{proof}

Next we prove that $|T_S|$ is small.

\begin{lemma} \label{SmallSTransversal}
	$|T_S| \leq 3m^{4/3}$.
\end{lemma}
\begin{proof}
	Observe that $|T_S|$ is at most the number of deactivated or neutralized pairs in $\mathbb{P}$; we bound this latter number. 
	Every deactivated or neutralized pair for a light $k$-set $N_i$ has a unique anchor from $\Gamma_i$,   
	thus the number of these pairs is at most $m^{1/3}$.  Since Lemma \ref{boundLightKSets} proves that
	the number of light $k$-sets is at most $3m$, it follows that $|T_S| \leq 3m^{4/3}$.
\end{proof}

To prove $|T_L| = O(m^{5/3})$, we first develop the next tool which finally
applies the leverage the $3$-crosses provide via Lemma \ref{PPinsideThreeCross}. 

\begin{lemma} \label{AtMostOneArcLabeled-g}
	For any $0 \leq x \leq \ell$, any $g \in G$, and any $y$ satisfying $0 \le y \leq t_x$, 
	there is at most one problematic pair with anchor $g$ and non-anchor $a_x^{(y)}$.
\end{lemma}
\begin{proof}  
	Suppose, to the contrary, that there are two 
	problematic private pairs,
	$p_{i_1}^{(j_1)}=\{g,a_x^{(y)}\}$ and $p_{i_2}^{(j_2)}=\{g,a_x^{(y)}\}$.
	This means, by definition, that $j_1,j_2 \leq y < t_x$.
	Because each $k$-set has at most one private pair with anchor $g$, it follows that $i_1 \neq i_2$.
	Without loss of generality, $i_1 < i_2$.  This means $t_{i_1} \geq t_{i_2}$,
	so $g$ is not a private vertex of $N_{i_2}^{(j_2)}$.
	Since $p_{i_1}^{(j_1)}$ is private for $N_{i_1}^{(j_1)}$, it follows that $t_{i_2} < j_1 \leq t_{i_1}$.
	Because $j_1 \leq y < t_x$,  we conclude that $y \not\in \{t_{i_2},t_x\}$.
      
    Now apply Lemma \ref{PPinsideThreeCross} to $p_{i_2}^{(j_2)}=\{g,a_x^{(y)}\}$ to obtain a $3$-cross 
    with  respect to $G_{i_2}^{(j_2-1)}$
    containing $p_{i_2}^{(j_2)}$.  Recall that this $3$-cross $C$ satisfies 
    $$p_{i_2}^{(j_2)} \subset C \subseteq
    \left( N_{i_2} \cup \{a_{i_2}^{(0)},\ldots,a_{i_2}^{(j_2)}\} \right) \setminus G_{i_2}^{(j_2-1)}$$
    \noindent
    Because $C$ is a $3$-cross, it must satisfy $C \cap \overline{N_{i_2}} \neq \varnothing$.  However,
    $p_{i_2}^{(j_2)} \subset N_{i_2}$ and $C$ contains $p_{i_2}^{(j_2)}$.  Consequently, we conclude that
    $C = \{a_{i_2}^{(\alpha)},g,a_x^{(y)}\}$, for some $a_{i_2}^{(\alpha)} \in \{a_{i_2}^{(0)},\ldots,a_{i_2}^{(j_2)}\}$.
    Since $t_{i_1} > t_{i_2}$, it follows that $a_{i_2}^{(\alpha)} \in N_{i_1}$.  Therefore, $C \subset  N_{i_1}$, contradicting
    that it is a $3$-cross (which requires  $C \cap \overline{N_{i_1}} \neq \varnothing$).
\end{proof}

\begin{lemma} \label{LightSetsHaveSmallTransversal}
	$|T_L| \leq 3 m^{5/3}$.
\end{lemma}
\begin{proof}
	Consider now an arbitrary light $k$-set $N_x$.
	Let $a_x^{(y_1)},\ldots,a_x^{(y_r)}$, with $y_1<\cdots<y_r$, be the collection of all kernel vertices in $\overline{N_x}$ that
	are the head of an arc of $D$ whose tail is in $T_L$.  We seek to give an upper bound on the number of tails of these arcs since
	this gives an upper bound for the number of vertices in $T_L$ contributed by $N_x$, due to Lemma \ref{Trichotomy}\ref{DifficultCase} and Lemma \ref{AtMostOneArcLabeled-g}.

	Suppose, for some $1 \leq w \leq r$, that $a_i^{(j)} \stackrel{g}{\longrightarrow} a_x^{(y_w)}$ is an arc of $D$ and $a_i^{(j)} \in T_L$.
	As noted right after the definition of $T_L$, 
	this means that $|I_g| = 1$.   Lemma \ref{Trichotomy}\ref{DifficultCase} 
	gives $I_g = \{i\}$ or $I_g = \{x\}$.
	If $I_g = \{i\}$, then $p_i^{(j)} = \{g,a_x^{(y_w)}\}$ is a troublesome pair.  Since we are assuming that $a_i^{(j)} \in T_L$, it follows that
	$a_i^{(j)} \not\in T_S$ so $p_i^{(j)}$ is an advanced pair.  However, this implies that 
	$y_w = t_x$ which contradicts condition \ref{IndexArcProperty} that
	$a_i^{(j)} \stackrel{g}{\rightarrow} a_x^{(y_w)}$  is  an arc of $D$.  So $I_g = \{x\}$.

	In other words, if $N_x$ is a light $k$-set, and  for $w$ fixed, $a_i^{(j)} \stackrel{g}{\longrightarrow} a_x^{(y_w)}$ is an arc of $D$ with $a_i^{(j)} \in T_L$, then $g \in \Gamma_x$. 	
	Lemma \ref{AtMostOneArcLabeled-g} states that there is at most one arc
	entering $a_x^{(y_w)}$ with label $g$.    Since $N_x$ is light, $|\Gamma_x| \leq m^{1/3}$, thus at most \
	$|\Gamma_x| \cdot r\leq  m^{1/3}\cdot r$ arcs with tail in $T_L$ enter $a_x^{(y_w)}$, for $w=1,\ldots, r$.	
	If $r \leq m^{1/3}$, then the number of vertices in $T_L$ contributed by $N_x$ would be at most $m^{2/3}$.
	Therefore, because Lemma \ref{boundLightKSets} shows $q \leq 3m$, to prove this lemma's conclusion ($|T_L| \leq 3 m^{5/3}$) it suffices to prove $r \leq m^{1/3}$.

	Toward this latter end, define, for $1 \leq w \leq r$,
	$$\Delta_x(w) = \{ g : \hbox{$a_i^{(j)} \stackrel{g}{\longrightarrow} a_x^{(y_w)}$ is an arc of $D$ and $a_i^{(j)} \in T_L$} \}.
	$$
	By the definition of $a_x^{(y_1)},\ldots,a_x^{(y_r)}$, we have $\Delta_x(w) \neq \varnothing$, for $w =1,\ldots, r$. 
	The prior paragraphs shows $\Delta_x(w) \subseteq \Gamma_x$, for all $1 \leq w \leq r$.  
	Consider an arbitrary $g\in \Delta_x(w)$.  By definition this means $D$ contains an arc
	$a_i^{(j)} \stackrel{g}{\longrightarrow} a_x^{(y_w)}$ and $a_i^{(j)} \in T_L$.
	Note that \ref{LightKSetException} guarantees that $g \not\in \Lambda_x^{(y_w)}$; hence $g \not\in G_x^{(y_w-1)}$.
	Therefore, $\Delta_x(w) \subseteq \Gamma_x \setminus  G_x^{(y_w-1)}$.
	
	Because 
	 $|\Gamma_x| \leq m^{1/3}$, to prove $r \leq m^{1/3}$ it is enough to prove:
	\begin{eqnarray}
		\Gamma_x \setminus  G_x^{(y_1-1)} \supsetneq \Gamma_x \setminus  G_x^{(y_2-1)} \supsetneq \cdots \supsetneq \Gamma_x \setminus  G_x^{(y_r-1)}.
	\end{eqnarray}
	
	Consider again an arbitrary $g\in \Delta_x(w) \subseteq \Gamma_x \setminus G_x^{(y_w-1)}$, for some $w \in \{1,\ldots,r-1\}$,
	with arc
	$a_i^{(j)} \stackrel{g}{\longrightarrow} a_x^{(y_w)}$ in $D$.
	If $g$ is not private for $N_x$ at time $y_w$, say $g \in N_u^{(y_w)}$ with $u \neq x$, then
	$a_x^{(y_w)} \in N_u$ contradicting that $p_i^{(j)} = \{g, a_x^{(y_w)}\}$ is a private pair for $N_i^{(j)}$.
	So $g$ is a private vertex for $N_x$ at time $y_w$ and  $g \not\in G_x^{(y_w-1)}$.  In other words,
	private vertices are available to produce the private pair $p_{x}^{(y_w)}$.
	Because private pair selection prioritizes vertices from $\Gamma_x \setminus G_x^{(y_w-1)}$ that are private to $N_x$ at time $y_w$,
	some vertex in $\Gamma_x \setminus G_x^{(y_w-1)}$ was chosen as the anchor for $p_{x}^{(y_w)}$.  This private vertex 
	is absent from $\Gamma_x \setminus G_x^{(y_w)}$ since it is added to $G_x^{(y_w-1)}$ to produce $G_x^{(y_w)}$.  
	Because $y_{w+1}-1 \geq y_w$,
	we conclude,
		$\Gamma_x \setminus  G_x^{(y_w-1)} \supsetneq \Gamma_x \setminus G_x^{(y_{w+1}-1)}$.
\end{proof}

\begin{theorem} \label{BoundOnT}
	$|T| \leq 6 m^{5/3} + 3m^{4/3} + 6m$.
\end{theorem}
\begin{proof}  Simply observe that by definition $T =  T_W \cup T_H \cup T_S \cup T_L$.  Now apply
	$|T_W| \leq 6m$ (Lemma \ref{SmallWTransversal}), $|T_H| \leq 3 m^{5/3}$ (Lemma \ref{SmallHTransversal}), $|T_S| \leq 3m^{4/3}$ (Lemma \ref{SmallSTransversal}), and
	$|T_L| \leq 3 m^{5/3}$ (Lemma \ref{LightSetsHaveSmallTransversal}).
\end{proof}

\section{A skew cross-intersecting system} 
\label{Skew}

In this section we apply the following theorem, first proven by Frankl \cite{Frankl} (see also \cite{Kalai}); it is the skew version of a theorem due to Bollob\'as \cite{Bollobas}.  This theorem is also 
presented in the book by Babai and Frankl (\cite{BabaiFrankl}, pages $94$--$95$).

\begin{theorem} \label{SkewBollobas} (Bollob\'as's Theorem - Skew Version)
If $A_1,\ldots,A_h$ are $r$-element sets and $B_1,\ldots,B_h$ are $s$-element sets such that
\begin{enumerate}[label=(\alph{enumi})]
\item \label{BollobasConditionA} $A_i$ and $B_i$ are disjoint for $i=1,\ldots,h$,
\item \label{BollobasConditionB} $A_i$ and $B_j$ intersect whenever $1 \leq i < j \leq h$
\end{enumerate}
then $h \leq {r+s \choose r}$.
\end{theorem}

A linearly ordered collection of pairs of sets, $\{(A_i,B_i)\}_{i=1}^h$, satisfying the hypotheses of Theorem \ref{SkewBollobas} is called a {\em skew intersecting set pair $(r,s)$-system}; abbreviate this to {\em skew $(r,s)$-system}.

Theorem \ref{SkewBollobas}  will be applied here  with $r=2$ and $ s=m$ 
 to obtain a skew $(2,m)$-system, where the $2$-sets are  members
 in the set $\mathbb{P}$  of  free pairs 
specified by the set $F\subset V\setminus T$ in Section \ref{transversal}, and  the  $m$-sets are derived iteratively from  the $k$-sets  corresponding to all the pairs in $\mathbb{P}$ as follows. 

Recall that  each $k$-set has a private pair at each stage until the $k$-set survives. First to every $N_i$ associate $t_i+1$ $m$-sets denoted $M_i^{(0)},\ldots,M_i^{(t_i)}$.
At stage $0$, set $M_i^{(0)} = \overline{N_i}$, for all $i=1,\ldots,\ell_0$.
For $i=1,\ldots,\ell$ and $j = 1,\ldots,t_i$, recursively define 
\begin{eqnarray} \label{RecursiveDefinitionOfMs}
	M_i^{(j)} = \left\{ 
	\begin{array}{lll}
		\left(M_i^{(j-1)} \setminus \{a_i^{(j-1)}\}\right) \cup \left\{g_i^{(j-1)}\right\} 
		                & \quad &\mbox{if } p_i^{(j-1)}= \left\{g_i^{(j-1)},u_i^{(j-1)}\right\}\mbox{ is free}\\ \\ 
		M_i^{(j-1)} & \quad & \mbox{if } p_i^{(j-1)} \mbox{ is not free}
	\end{array}
	\right.
\end{eqnarray}

Note that, because  $a_i^{(j-1)}\in M_i^{(j-1)}$, $g_i^{(j-1)}\not\in M_i^{(j-1)}$, and $|M_i^{(0)}|=m$, it follows that $|M_i^{(j)}| = m$, for all 
$i=1,\ldots,\ell$ and $j = 1,\ldots,t_i$.  
This
recursive process will never remove $a_i^{(t_i)}$ from $M_i^{(0)}$
because the process halts at stage $j=t_i$.

Now define the set-pair system
$${\cal F} = \left\{ (p_i^{(j)},M_i^{(j)}): p_i^{(j)} \in \mathbb{P} \mbox{ is free}\right\},$$
where ${\cal F}$ is ordered linearly and chronologically via lexicographical order: 
$$(p_i^{(j)},M_i^{(j)}) < (p_x^{(y)},M_x^{(y)}) \qquad \Longleftrightarrow \qquad \mbox{($j < y$) or ($j=y$ and $i < x$}).$$

\begin{theorem} \label{SkewSystem} ${\cal F}$ is a skew $(2,m)$-system.
\end{theorem}
\begin{proof} Clearly $|p_i^{(j)}| = 2$ and $|M_i^{(j)}| = m$, for
all $(p_i^{(j)},M_i^{(j)}) \in {\cal F}$.  Because $p_i^{(j)}=\{ g_i^{(j)},u_i^{(j)}\}$ is private to $N_i$ at time $j$,
it follows that $p_i^{(j)} \subset N_i$; so, $p_i^{(j)} \cap M_i^{(0)} = \varnothing$.  Observe 
that $g_i^{(j)}\notin M_i^{(j)}$,  
because $g_i^{(j)}$ is added to $M_i^{(j)}$ at time $j+1$ 
by the recursive process generating the $m$-sets.  
Moreover,  
$\{g_i^{(j)},u_i^{(j)}\}\cap \{g_i^{(0)},\ldots, g_i^{(j-1)}\}=\varnothing$, by the iterative choice of private pairs, 
thus $u_i^{(j)}\notin M_i^{(0)}$ implying that $u_i^{(j)}\notin M_i^{(j)}$. 
Therefore, $p_i^{(j)} \cap M_i^{(j)} = \varnothing$, showing that hypothesis \ref{BollobasConditionA} is satisfied in Theorem \ref{SkewBollobas}.

Now we prove the system satisfies hypothesis \ref{BollobasConditionB}. 
Suppose $(p_i^{(j)},M_i^{(j)}), (p_x^{(y)},M_x^{(y)}) \in {\cal F}$ and 
$(p_i^{(j)},M_i^{(j)}) < (p_x^{(y)},M_x^{(y)})$, in particular, $j\leq y$.  
We must prove $p_i^{(j)} \cap M_x^{(y)} \neq \varnothing$.

If $i= x$, then $j < y$ so $g_i^{(j)}\in M_i^{(j+1)} = M_x^{(j+1)}$ because $p_i^{(j)}$ is free
and therefore (\ref{RecursiveDefinitionOfMs}) swaps $g_i^{(j)}$ into $M_i^{(j+1)}$.  In this case, $g_i^{(j)}\in M_x^{(y)}$. 
Consequently, we may assume $i \neq x$.

If $g_i^{(j)} \in M_x^{(y)}$, then $p_i^{(j)} \cap M_x^{(y)} \neq \varnothing$.   
So we may assume
$g_i^{(j)} \notin M_x^{(y)}$.  
Elements from $G$ are
only added to $M_x^{(0)}$ to get to $M_x^{(y)}$, so $g_i^{(j)} \notin M_x^{(0)} = \overline{N_x}$ which implies $g_i^{(j)} \in N_x$.  
Since $j\leq y$ and $p_i^{(j)}$ is private to $N_i$ at time $j$, we conclude that 
$u_i^{(j)} \notin N_x$ because $N_x$ survives at stage $j$ (that is, $t_x \geq y \geq j$).  
So $u_i^{(j)} \in \overline{N_x} = M_x^{(0)}$. 
   
 If $u_i^{(j)} \in G$, or $u_i^{(j)} = a_x^{(s)}$ and $s \geq y$,
then $a_x^{(s)}$ is not removed from $M_x^{(0)}$ during the process generating $M_x^{(y)}$,
meaning $p_i^{(j)} \cap M_x^{(y)} \neq \varnothing$.
So we may assume that 
$u_i^{(j)} = a_x^{(s)}$, for some $s  \leq y-1$.  The private pair selection and the promotion process (\ref{prom}) guarantees that $j \leq s$.

Setting $g = g_i^{(j)}$, we find 
$p_i^{(j)} = \{g, a_x^{(s)}\}$ with $g \in N_i \cap (N_x \cap G)$ and $j \leq s \leq y - 1 < y \leq t_x$.
In other words, according to Definition \ref{problematic},  $p_i^{(j)}$ is a problematic pair.
Lemma \ref{Trichotomy}\ref{Tri} yields $|I_g| = 1$ or $|I_g|=2$.
If $|I_g|=2$, then $a_i^{(j)} \in T_W$ which contradicts the assumption that $p_i^{(j)}$ is free.
So we may assume $|I_g| = 1$.  Lemma \ref{CriticalIndexSet}\ref{CriticalIndexSetsOfArcs} implies
that $I_g = \{i\}$ or  $I_g = \{x\}$.

If $N_x$ is heavy, then the definition of $T_H$ means $a_x^{(s)} \in T_H$, so $p_x^{(s)}$ is not free.  
Accordingly, $a_x^{(s)}$ is never removed from $M_x^{(0)}$ in the production of $M_x^{(y)}$.  Hence, 
$a_x^{(s)}  \in p_i^{(j)} \cap M_x^{(y)}$.  So we may assume that $N_x$ is light.

If $g \in \Lambda_x^{(y)}$, 
then $g$ is the anchor for
some non-neutralized private pair $p_x^{(\beta)}$, for some $0 \leq \beta \leq y-1$ and $g \in \Gamma_x$.
In this case, Property~\eqref{property:L} guarantees that $a_x^{(\beta)} \not\in T$.
It follows that $p_x^{(\beta)}$ is free
so (\ref{RecursiveDefinitionOfMs}) swaps $g$ into $M_x^{(\beta+1)}$.  This means $g \in M_x^{(y)}$ yielding
$g  \in p_i^{(j)} \cap M_x^{(y)}$.   
So we may assume that $g \not\in \Lambda_x^{(y)}$.

Because $p_i^{(j)}$ is free, we conclude that $a_i^{(j)} \not\in T_S$ so $p_i^{(j)}$ is active.

To summarize, we now have these remaining conditions: $p_i^{(j)} = \{g, a_x^{(s)}\}$, $p_i^{(j)}$ is active,
$g \in N_i \cap N_x$, $j \leq s < y  \leq t_x$, $N_x$ is light, and $g \not\in \Lambda_x^{(y)}$.
These conditions guarantee 
$a_i^{(j)} \stackrel{g}{\longrightarrow} a_x^{(s)}$  is an arc of $D$. 
But this case can not occur because
the definition of $T_L$ would give $a_i^{(j)} \in T_L$, implying $a_i^{(j)} \not\in F$.
This contradicts the assumption that $p_i^{(j)}$ is free.
\end{proof}


Now we may state the main theorem of the paper, a new upper bound on the order of an \nmsystem{}
\begin{theorem} \label{maintheorem}
Any \nmsystem{} satisfies $$n\leq  {m+2 \choose 2} + 6 m^{5/3} + 3m^{4/3} + 9m - 3.$$
\end{theorem}
\begin{proof} 
Recall that ${F} = A\setminus T$ is a free set of pairs; so
$ |{F}| = |A| - |T|$, and 
the skew Bollob\'{a}s theorem yields
$
|{F}| \leq {m+2 \choose 2}.
$   
Therefore, $|A| \leq {m+2 \choose 2} + |T|$. 
By Lemma \ref{FirstLemma}(\ref{GUpperBound}), $|G|\leq 3(m-1)$; and using the bound on $|T| \leq 6 m^{5/3} + 3m^{4/3} + 6m$ in Theorem \ref{BoundOnT} we obtain
\begin{eqnarray*} 
n & = & |G| + |A| \\
  &\leq & |G| + {m+2 \choose 2} + |T| \\ 
	& \leq & 3(m-1) + {m+2 \choose 2} + 6 m^{5/3} + 3m^{4/3} + 6m \\
	& \leq & {m+2 \choose 2} + 6 m^{5/3} + 3m^{4/3} + 9m - 3.
\end{eqnarray*}
\end{proof}


\section*{Acknowledgments}
We thank Adam Jobson  and Tim Pervenecki for their helpful remarks and discussions during the project. We are also very grateful to Zsolt Tuza for his valuable comments on an earlier version of this paper and bringing our attention to reference \cite{Tuza}.

\end{document}